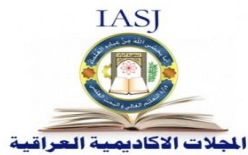

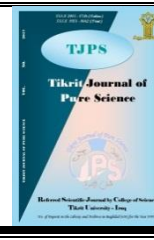

# On Harmonic Univalent Functions Defined by Dziok-Srivastava Operator

**Mays S. Abdul Ameer[1], Abdul Rahman S. Juma[2], Raheam A. Al-Saphory[1]**

*[1]Department of Mathematics, College of Education for Pure Science, Tikrit University, Tikrit, Iraq*
*[2] Department of Mathematics, College of Education for Pure Science, University of Anbar, Ramadi, Iraq*

https://doi.org/10.25130/tjps.v26i1.110

## ABSTRACT

The purpose of this work is to present a class of harmonic univalent functions defined by the Dziok-Srivastava operator. Some geometric properties like coefficients conditions, distortion theorem, convolution (Hadamard product), convex combination and extreme points are investigated.

## 1. Introduction

Let A denote the class of functions $f(z)$ which are analytic in the open unit disk
$U = \{z: |z| < 1\}$. Each $f \in A$ can be represented by $f = h + \overline{g}$, with $h$ and $g$ of analytic type in U. We say that $h$ is an analytic part and $g$ the related co-analytic part of $f$ (see[2]). Thus for
$f = h + \overline{g} \in A$, we can write
$h(z) = z + \sum_{n=2}^{\infty} a_n z^n$, $g(z) = \sum_{n=1}^{\infty} b_n z^n$,
$(0 \leq |b_1| < 1)$ ...... (1)
Let $f = h + \overline{g}$ given by (1) and $F_{(\lambda \rho)(\mu q),b}^{s,a,\lambda}$ is Dziok-Srivastava operator of $f$ and is given by [8]
$F_{(\lambda \rho)(\mu q),b}^{s,a,\lambda} f(z) =$

$z + \sum_{n=2}^{\infty} \frac{\prod_{j=1}^{p}(\lambda_j + 1)_{n-1}}{\prod_{j=1}^{q}(\mu_j + 1)_{n-1}} (\frac{a+1}{a+n})^s \frac{\Delta(a+n,b,s,\lambda)}{\Delta(a+1,b,s,\lambda)} a_n \frac{z^n}{n!}$

.... (2)
Where

$(p, q \in N_O;$ $\lambda_j \in C$, $(j = 1,2, \dots, p);$ $a, \mu_j \in \frac{C}{Z_O^-}$,
$(j = 1,2, \dots, q); Z_O^- = \{0,1,2 \dots\};$
$s, z \in C)$,
but
$\Phi_n(s, a; b, \lambda) = \frac{\prod_{j=1}^{p}(\lambda_j + 1)_{n-1}}{\prod_{j=1}^{q}(\mu_j + 1)_{n-1} n!} (\frac{a+1}{a+n})^s \frac{\Delta(a+n,b,s,\lambda)}{\Delta(a+1,b,s,\lambda)}.$
...... (3)

Let $B_H(\lambda, \alpha, k)$ be the family of functions harmonic type of form (1) by

$Re \left\{ \frac{F_{(\lambda \rho)(\mu q),b}^{s,a,\lambda} f(z)}{z(F_{(\lambda \rho)(\mu q),b}^{s,a,\lambda} f(z))'} \right\} > \alpha$, $0 \leq \alpha < 1$,

where $F_{(\lambda \rho)(\mu q),b}^{s,a,\lambda} f(z)$ is defined by (2). .... (4)
Let $\overline{B_H}(\lambda, \alpha, k)$ denoted that the subclass of $B_H(\lambda, \alpha, k)$ such that $h$ and $g$ are the from
$h(z) = z - \sum_{n=2}^{\infty} a_n z^n$, $g(z) = (-1)^n \sum_{n=1}^{\infty} b_n z^n, |b_1| < 1.$ ... (5)

## 2. The Main Results

In this section, the main important results are stated and proved an enough coefficient to functions of harmonic univalent types.

**Theorem 2.1:** Let $f = h + \overline{g}$ be given by (1), if
$\sum_{n=2}^{\infty} (1 - \alpha n)|\Phi_n(s, a; b, \lambda)||a_n||z|^n + \sum_{n=1}^{\infty} (1 - \alpha n)|\Phi_n(s, a; b, \lambda)||b_n||\overline{z}|^n \leq (1 - \alpha)$, ...... (6)
where $a_1 = 1, 0 \leq \alpha < 1$, then, $f$ harmonic of sense-preserving type to $U$ with
$\Phi_n(s, a; b, \lambda) = \frac{\prod_{j=1}^{p}(\lambda_j + 1)_{n-1}}{\prod_{j=1}^{q}(\mu_j + 1)_{n-1} n!} (\frac{a+1}{a+n})^s \frac{\Delta(a+n,b,s,\lambda)}{\Delta(a+1,b,s,\lambda)}.$

**Proof:** For $z_1 \neq z_2$, we have
$\left| \frac{f(z_1) - f(z_2)}{h(z_1) - h(z_2)} \right| \geq 1 - \left| \frac{g(z_1) - g(z_2)}{h(z_1) - h(z_2)} \right|$

$= 1 - \left| \frac{\sum_{n=1}^{\infty} b_n(z_1^n - z_2^n)}{(z_1 - z_2)\sum_{n=2}^{\infty} a_n(z_1^n - z_2^n)} \right| \geq 1 - \frac{\sum_{n=1}^{\infty} n|b_n|}{1 - \sum_{n=2}^{\infty} n|a_n|}$

$\geq 1 - \frac{\frac{1-\alpha}{(1-\alpha n)|\varphi_n(s,a;b,\lambda)|}|b_n|}{\frac{1-\alpha}{(1-\alpha n)|\varphi_n(s,a;b,\lambda)|}|a_n|} \geq 0$ .

Hence proved univalent, since

$\left| h'(z) \right| \geq 1 - \sum_{n=2}^{\infty} n|a_n||z|^{n-1}$

$> 1 - \sum_{n=2}^{\infty} n|a_n|$

$\sum_{n=2}^{\infty} \frac{(1-\alpha n)|\varphi_n(s,a;b,\lambda)|}{1-\alpha}|a_n|$

$\geq 1 - \sum_{n=1}^{\infty} \frac{(1-\alpha n)|\varphi_n(s,a;b,\lambda)|}{1-\alpha}|b_n| \geq$

$\sum_{n=1}^{\infty} n|b_n|$

$> \sum_{n=1}^{\infty} n|b_n||z|^{n-1} \geq |g'(z)|$ .

By using the case of $Re(w) \geq \alpha, \Leftrightarrow |1 - \alpha + w| \geq |1 + \alpha - w|$

for $(0 \leq \alpha < 1)$, it show that

$|A(z) + (1-\alpha)B(z)| - |A(z) - (1+\alpha)B(z)| \geq 0$, ….(7)

where

$A(z) = F_{(\lambda \rho)(\mu q),b}^{s,a,\lambda} f(z)$ , $B(z) = z(F_{(\lambda \rho)(\mu q),b}^{s,a,\lambda} f(z))$.

Consider $|A(z) + (1-\alpha)B(z)|$

$\leq \sum_{n=2}^{\infty} (1 + n - \alpha n)|\varphi_n(s,a;b,\lambda)||a_n||z|^n + \sum_{n=1}^{\infty} (1 + n - \alpha n)|\varphi_n(s,a;b,\lambda)||b_n||\bar{z}|^n$

$+2 - \alpha$ . ………(8)

And

$|A(z) - (1+\alpha)B(z)|$

$=$

$|z + \sum_{n=2}^{\infty} \Phi_n(s,a;b,\lambda)a_n z^n + (-1)^k \sum_{n=1}^{\infty} \Phi_n(s,a;b,\lambda)b_n \bar{z}^n - (1+\alpha)[z +$

$\sum_{n=2}^{\infty} n\Phi_n(s,a;b,\lambda)a_n z^n +$

$(-1)^k \sum_{n=1}^{\infty} n\Phi_n(s,a;b,\lambda)b_n \bar{z}^n|$

$\leq \sum_{n=2}^{\infty} (1 - n - \alpha n)|\Phi_n(s,a;b,\lambda)||a_n||z|^n + \sum_{n=1}^{\infty} (1 - n - \alpha n)|\Phi_n(s,a;b,\lambda)||b_n||\bar{z}|^n$

$+\alpha$. ……(9)

Substituting (8) and (9) in (7), we obtain

$2\sum_{n=1}^{\infty} (1 - \alpha n)|\Phi_n(s,a;b,\lambda)||a_n||z|^n +$

$2\sum_{n=1}^{\infty} (1 - \alpha n)|\Phi_n(s,a;b,\lambda)||b_n||\bar{z}|^n$

$+2(1 - \alpha)$,

so

$\sum_{n=2}^{\infty} (1 - \alpha n)|\Phi_n(s,a;b,\lambda)||a_n||z|^n + \sum_{n=1}^{\infty} (1 - \alpha n)|\Phi_n(s,a;b,\lambda)||b_n||\bar{z}|^n$

$\leq (1 - \alpha)$.

**Theorem 2.2:** Let $f = h + \overline{g}$ is formulated by (5). Then $f(z) \in \overline{B_H}(\lambda, \alpha, k)$ if and only if

$\sum_{n=2}^{\infty} (1 - \alpha n)|\Phi_n(s,a;b,\lambda)||a_n||z|^n + \sum_{n=1}^{\infty} (1 - \alpha n)|\Phi_n(s,a;b,\lambda)||b_n||\bar{z}|^n$

$\leq (1 - \alpha)$, ….(10)

where

$\varphi_n(s,a;b,\lambda) = \frac{\prod_{j=1}^{p}(\lambda_j + 1)_{n-1}}{\prod_{j=1}^{q}(\mu_j + 1)_{n-1}n!}\left(\frac{a+1}{a+n}\right)^s \frac{\Lambda(a+n,b,s,\lambda)}{\Lambda(a+1,b,s,\lambda)}$,

and $\lambda_j \epsilon C (j = 1, \dots, p), \mu_j \epsilon \frac{c}{z_0} (j = 1, \dots, q)$

**Proof:** Since $\overline{B_H}(\lambda, \alpha, k) \subset B_H(\lambda, \alpha, k)$ , we just need to prove the only if part of the theorem. We notice that the condition (5) is equation to

$Re\left\{ \frac{F_{(\lambda \rho)(\mu q),b}^{s,a,\lambda} f(z)}{z(F_{(\lambda \rho)(\mu q),b}^{s,a,\lambda} f(z))} \right\} > \alpha$ , $0 \leq \alpha < 1$

$Re\left[ \frac{z - \sum_{n=2}^{\infty} |\Phi_n(s,a;b,\lambda)||a_n|z^n + (-1)^{2s+1}\sum_{n=1}^{\infty} |\Phi_n(s,a;b,\lambda)||b_n|\bar{z}^n}{z - \sum_{n=2}^{\infty} n|\Phi_n(s,a;b,\lambda)||a_n|z^n + (-1)^{2s}\sum_{n=1}^{\infty} n|\Phi_n(s,a;b,\lambda)||b_n|\bar{z}^n} \right] > \alpha$.

$Re\left[ \frac{z - \sum_{n=2}^{\infty} |\Phi_n(s,a;b,\lambda)||a_n|z^n + (-1)^{2s+1}\sum_{n=1}^{\infty} |\Phi_n(s,a;b,\lambda)||b_n|\bar{z}^n}{z - \sum_{n=2}^{\infty} n|\Phi_n(s,a;b,\lambda)||a_n|z^n + (-1)^{2s}\sum_{n=1}^{\infty} n|\Phi_n(s,a;b,\lambda)||b_n|\bar{z}^n} \right] \geq 0$.

$Re\left[ \frac{(1-\alpha) - \sum_{n=2}^{\infty} (1-\alpha n)|\Phi_n(s,a;b,\lambda)||a_n|z^{n-1} - \frac{\bar{z}}{z}\sum_{n=1}^{\infty} (1-\alpha n)|\Phi_n(s,a;b,\lambda)||b_n|\bar{z}^{n-1}}{1 - \sum_{n=2}^{\infty} |\Phi_n(s,a;b,\lambda)||a_n|z^{n-1} + \frac{\bar{z}}{z}\sum_{n=1}^{\infty} |\Phi_n(s,a;b,\lambda)||b_n|\bar{z}^{n-1}} \right] \geq 0$. …(11)

The condition (11) must satisfy for all values of $z$ on $|z| \in (0, \mu)$, we must have the positive real axis, where

$Re\left[ \frac{(1-\alpha) - \sum_{n=2}^{\infty} (1-\alpha n)|\Phi_n(s,a;b,\lambda)||a_n|\mu^{n-1} - \frac{\bar{z}}{z}\sum_{n=1}^{\infty} (1-\alpha n)|\Phi_n(s,a;b,\lambda)||b_n|\bar{\mu}^{n-1}}{1 - \sum_{n=2}^{\infty} |\Phi_n(s,a;b,\lambda)||a_n|\mu^{n-1} + \frac{\bar{z}}{z}\sum_{n=1}^{\infty} |\Phi_n(s,a;b,\lambda)||b_n|\bar{\mu}^{n-1}} \right] \geq 0$. …(12)

If the condition (6) does not hold then the numerator in (8), when goes to 1 is negative. This is a contradiction with the situation case where $f \in \overline{B_H}(\lambda, \alpha, k)$ and so the proof is accomplished.

## 3. The Distortion Theorem

**Theorem 3.1:** Suppose $f \in \overline{B_H}(\alpha, \beta, \lambda)$. Then $|z| = r < 1$, we can get

$|f(z)| \leq (1 + |b_1|)r + \frac{1}{\Phi_2(s,a;b,\lambda)}\left[\frac{1-2\alpha}{1-\alpha} - \frac{1+2\alpha}{1-\alpha}|b_1|\right]r^2$ , $|z| = r < 1$

and

$|f(z)| \geq (1 + |b_1|)r - \frac{1}{\Phi_2(s,a;b,\lambda)}\left[\frac{1-2n}{1-\alpha} - \frac{1+2n}{1-\alpha}|b_1|\right]r^2$, $|z| = r < 1$.

**Proof:** We have

$|f(z)| \leq (1 + |b_1|)r + [(|a_n| + |b_n|)]r^n$

$\leq (1 + |b_1|)r + [(|a_n| + |b_n|)]r^2$

$= (1 + |b_1|)r + \frac{1-\alpha}{\Phi_2(s,a;b,\lambda)(1-2\alpha)}\left[\frac{1-2\alpha}{1-\alpha}|a_n| + \frac{1-2\alpha}{1-\alpha}|b_n|\right] |\Phi_2(s,a;b,\lambda)|r^2 \leq (1 + |b_1|)r + \frac{1}{|\Phi_2(s,a;b,\lambda)|(1-2\alpha)}\left(1 - \frac{1+2\alpha}{1-\alpha}|b_1|\right)r^2 \leq (1 + |b_1|)r + \frac{1}{|\Phi_2(s,a;b,\lambda)|}\left[\frac{1-\alpha n}{1-\alpha} - \frac{1+2\alpha}{1-\alpha}|b_1|\right]r^2$.

## 4. The Convolution (Hadamard product)

Let

$f(z) = z - \sum_{n=2}^{\infty} |a_n|z^n + (-1)^k \sum_{n=1}^{\infty} |b_n|\bar{z}^n$

and

$g(z) =$
$z - \sum_{n=2}^{\infty} \quad |c_n|z^n + (-1)^k \sum_{n=1}^{\infty} \quad |d_n|\overline{z}^n$.
Then, form the convolution of $f(z)$ and $g(z)$, we can obtain
$(f * g)(z) = f(z) * g(z) \qquad =$
$z - \sum_{n=2}^{\infty} \quad |a_n||c_n|z^n +$
$(-1)^k \sum_{n=1}^{\infty} \quad |b_n| |d_n|\overline{z}^n$.

**Theorem 4.1:** Let $f(z) \in \overline{B_H}(\lambda, \alpha, k)$ and $g(z) \in \overline{B_H}(\lambda, \beta, k)$. Then for
$0 \le \beta \le \alpha < 1$, we have
$(f * g)(z) \in \overline{B_H}(\lambda, \alpha, k) \subset \overline{B_H}(\lambda, \beta, k)$, then satisfy (6) and since ($|c_n| \le 1, |d_n| \le 1$),
we write
$\sum_{n=1}^{\infty} \quad (\frac{(1-\alpha n)}{1-\alpha}|a_n c_n| +$
$\frac{(1-\alpha n)}{1-\alpha}|b_n d_n|) |\Phi_n(s, a; b, \lambda)|$
$\le \sum_{n=1}^{\infty} \quad (\frac{(1-\alpha n)}{1-\alpha}|a_n| + \frac{(1-\alpha n)}{1-\alpha}|b_n|) |\Phi_n(s, a; b, \lambda)|$
The last inequality is bounded of the right hand side above by(1), then
$f * g \in \overline{B_H}(\lambda, \alpha, k) \subset \overline{B_H}(\lambda, \beta, k)$.

## 5. The Convex Combination
This section, is devoted to prove that the space $\overline{B_H}(\lambda, \alpha, k)$ is closed under convex combination. Suppose that the $f_i(z)$ is formulated for $i = 1, 2, 3, \ldots, m$, by the following form
$f_i(z) = z - \sum_{n=2}^{\infty} \quad |a_{n,i}|z^n + (-1)^k \sum_{n=1}^{\infty} \quad |b_{n,i}|\overline{z}^n$.
….(13)

**Theorem 5.1:** Assume that the functions $f_i$ is given by (13) be in the class $\overline{B_H}(\lambda, \alpha, k)$, for every $i = 1, 2, \ldots, m$. Then the functions $\tau_i(z)$ well-defined by
$\tau_i(z) = \sum_{n=1}^{\infty} \quad C_i f_{i(z)}, \ 0 \le c_i \le 1$ are also in the class $\overline{B_H}(\lambda, \alpha, k)$ where $\sum_{n=1}^{\infty} \quad C_i = 1$.

**Proof.** In view of the $t_i(z)$ definition, we can write reformulate
$\tau_i(z) = z - \sum_{n=2}^{\infty} \quad (\sum_{i=1}^{\infty} \quad C_i|a_{n,i}|)z^n +$
$(-1)^k \sum_{n=1}^{\infty} \quad (\sum_{i=1}^{\infty} \quad C_i|b_{n,i}|)\overline{z}^n$ …. (14)
Further, since $f_i(z)$ are in $\overline{B_H}(\lambda, \alpha, k)$, for every $i = 1, 2, \ldots, m$, then we have
$\sum_{n=2}^{\infty} \quad \frac{(1-\alpha)}{1-\alpha}(\sum_{i=1}^{\infty} \quad C_i|a_{n,i}|)z^n +$
$\sum_{n=1}^{\infty} \quad \frac{(1-\alpha n)}{1-\alpha}(\sum_{i=1}^{\infty} \quad C_i|b_{n,i}|)\overline{z}^n$
$= \sum_{i=1}^{\infty} \quad C_i(\sum_{n=2}^{\infty} \quad (1-\alpha n)|\varphi_n(s, a; b, \lambda)||a_{n,i}| + \sum_{n=1}^{\infty} \quad (1-\alpha n)|\varphi_n(s, a; b, \lambda)||b_{n,i}|)$
$\le \sum_{i=1}^{\infty} \quad C_i(1-\alpha) \le (1-\alpha)$.

## 6. The Extreme Point
In this part, we get the extreme points for the class $\overline{B_H}(\alpha, \beta, \lambda)$.

**Theorem 6.1:** Suppose $f$ be given by (5). Then $f \in \overline{B_H}(\alpha, \beta, \lambda), \Leftrightarrow$
$f(z) = \sum_{n=1}^{\infty} \quad (T_n h_n(z) + S_n g_n(z))$, …(15)
where
$h_1(z) = z$,
$h_n(z) = z - (\frac{(1-\alpha)}{(1-\alpha n)|\Phi_n(s; a; b, \lambda)|})z^n, \ n = 2, 3, \ldots$
and

$g_n(z) = z + (-1)^{k-1} (\frac{(1-\alpha)}{(1-\alpha n)|\Phi_n(s; a; b, \lambda)|})\overline{z}^n, \ n = 1, 2, \ldots$
$\sum_{i=1}^{\infty} \quad (T_n + S_n) = 1, \ T_n \ge 0$ and $S_n \ge 0$.
In particular, the extreme points of $f \in \overline{B_H}(\alpha, \beta, \lambda)$ are $\{h_n\}$ and $\{g_n\}$.
**Proof:** The form (15), we get
$f(z) = \sum_{n=1}^{\infty} \quad (T_n h_n(z) + S_n g_n(z))$
$=$
$\sum_{n=1}^{\infty} \quad (T_n + S_n)z - \sum_{n=2}^{\infty} \quad \frac{(1-\alpha)}{1-\alpha n|\Phi_n(s; a; b, \lambda)|}T_n z^n$
$\quad +(-1)^{n-1} \sum_{n=1}^{\infty} \quad \frac{((1-\alpha))}{1-\alpha n|\Phi_n(s; a; b, \lambda)|}S_n \overline{z}^n$
$= z - \sum_{n=2}^{\infty} \quad (\frac{(1-\alpha)}{(1-\alpha n)|\Phi_n(s; a; b, \lambda)|})T_n z^n$
$+(-1)^k \sum_{n=1}^{\infty} \quad (\frac{(1-\alpha)}{(1-\alpha n)|\Phi_n(s; a; b, \lambda)|})S_n \overline{z}^n$.
Therefore
$\sum_{n=2}^{\infty} \quad \frac{(1-\alpha n)|\Phi_n(s; a; b, \lambda)|}{1-\alpha}|a_n|$
$+\sum_{n=1}^{\infty} \quad \frac{(1-\alpha n)|\Phi_n(s; a; b, \lambda)|}{1-\alpha}|b_n|$
$\sum_{n=2}^{\infty} \quad \frac{(1-\alpha n)|\Phi_n(s; a; b, \lambda)|}{1-\alpha}(\frac{(1-\alpha)}{1-\alpha n|\Phi_n(s; a; b, \lambda)|} a_n) +$
$\sum_{n=1}^{\infty} \quad \frac{(1-\alpha n)|\Phi_n(s; a; b, \lambda)|}{1-\alpha}(\frac{(1-\alpha)}{1-\alpha n|\Phi_n(s; a; b, \lambda)|} b_n)$
$= \sum_{n=1}^{\infty} \quad T_n + \sum_{i=1}^{\infty} \quad S_n = 1 - T_n \le 1$.
And so $f \in \overline{B_H}(\alpha, \beta, \lambda)$.
Conversely, assume that $f \in \overline{B_H}(\alpha, \beta, \lambda)$.
Letting
$T_1 = 1 - \sum_{n=2}^{\infty} \quad T_n + \sum_{n=1}^{\infty} \quad S_n$.
Putting
$T_n = \frac{(1-\alpha n)|\Phi_n(s; a; b, \lambda)|}{1-\alpha}|a_n|, \ n = 2, 3, \ldots$
$S_n = \frac{(1-\alpha n)|\Phi_n(s; a; b, \lambda)|}{1-\alpha}|b_n|, \ n = 1, 2, \ldots$
we obtain
$f(z) = z - \sum_{n=2}^{\infty} \quad a_n z^n + (-1)^n \sum_{n=1}^{\infty} b_n \ \overline{z}^n$
$= z - \sum_{n=2}^{\infty} \quad \frac{(1-\alpha n)|\Phi_n(s; a; b, \lambda)|}{1-\alpha}T_n z^n +$
$(-1)^n \sum_{n=1}^{\infty} \quad \frac{(1-\alpha n)|\Phi_n(s; a; b, \lambda)|}{1-\alpha}S_n \overline{z}^n$.
$= z - \sum_{n=2}^{\infty} \quad [z - h_n(z)]T_n - \sum_{n=1}^{\infty} \quad [z - g_n(z)]S_n$
$=$
$1 - [\sum_{n=2}^{\infty} \quad T_n - \sum_{n=1}^{\infty} \quad S_n]z + \sum_{n=2}^{\infty} \quad T_n h_n(z) + \sum_{n=1}^{\infty} \quad S_n g_n(z)$
$= \sum_{i=1}^{\infty} \quad (T_n h_n(z) + S_n g_n(z))$.
This completes the proof.

## Conclusion
We have shown that a new class to functions of harmonic univalent type, interesting results concerning the harmonic univalent functions defined by the Dziok-Srivastava operator. Thus, some geometric properties like coefficients conditions, distortion theorem, convolution (Hadamard product), extreme points and convex combination are investigated and examined. Finally, Moreover, many problems still opened, for example, the extension of these results to the case of subclasses for various linear operator [9-11].

**References**

[1] Clunie J., Sheil T. (1984). Harmonic univalent functions, Ann. Acad. Sci. Fenn. Ser. A. I. Math.,vol. 9, pp 3-25.

[2] Cotirla L. I. (2009). Harmonic univalent functions defined by an integral operator, A. Univ. Apul., pp95-105.

[3] Gregorczyk M., Widomski J. (2010). Harmonic mappings in the exterior of the unit disk, annales universitatis mariae curie-sklodowska, Vol. lxiv, no. 1, pp63-73.

[4] Janangiri J.M. (1999). Harmonic functions starlike in the unit disc, J. Math. Anal. Appl.,
 Vol. 235, pp470–447.

[5] Juma A. R. S. (2010). On harmonic univalent function defined by generalized Salagean derivatives, acta, no.2 3, ,179-188.

[6] Juma A. R. S. and. Kulkarni S. R. (2007). On univalent functions with negative coefficients by using generalized Salagean operator, Filomat, Vol. 21,173-184.

[7] Srivastava H. M. and Gaboury S. (2015). A new class of analytic functions defined by means of a generalization of the Srivastava-Attiya Operator, J. Inequal. Appl. 2015, Article ID 39,
pp1–15.

[8] Srivastava H. M., Juma A. R. S. and Zayed H. M. (2018). On Univalence conditions for an integral operator defined by a generalization of the Srivastava-Attiya operator. Filomat, vol. 32. no. 6, pp2101–2114.

[9] Mahmoud M., Juma A. R. and Al-Saphory R., On bi-univalent functions involving Srivastava-Attiya operator, Italian Journal of Pure and Applied Mathematics, Accepted 2020.

[10] Abdul Ameer M., Juma A. R. and Al-Saphory R., Differential Subordination for Higher-Order Derivatives of Multivalent Functions, Journal of Physics: Conference Series, Accepted 2020.

[11] Abdul Ameer M., Juma A. R. and Al-Saphory R., Harmonic Multivalent Functions Defined by General Integral Operator, Al-Qadisiyah Journal of Pure Science, vol. 26, no. 1, ppMath 1-12, 2021.

## الدوال الاحادية التكافؤ التوافقية معرفة بواسطة المؤثر Dziok-Srivastava

ميس صالح عبد الامير [1]، عبد الرحمن سلمان جمعة [2] ، رحيم احمد منصور [1]

[1] قسم الرياضيات ،كلية التربية للعلوم الصرفة ، جامعة تكريت ، تكريت ، العراق

[2] قسم الرياضيات ،كلية التربية للعلوم الصرفة ، جامعة الانبار ، الرمادي ، العراق

**الملخص**

الغرض من هذا العمل هو تقديم فئة من الدوال التوافقية أحادية التكافؤ التي حددها عامل التشغيل Dziok-Srivastava.

تم دراسة بعض الخصائص الهندسية مثل شروط المعاملات، نظرية التشويه، الالتواء (ضرب هادمارد)، التركيبة المحدبة والنقاط المتطرفة.

.